\documentclass[a4paper,11pt,twoside,reqno]{amsart}

\usepackage[utf8]{inputenc}
\usepackage[plainpages=false,pdfpagelabels=true]{hyperref}
\usepackage{amssymb,amsthm}
\usepackage[margin=1in]{geometry}
\usepackage{slashed}
\usepackage{comment}
\usepackage{caption}

\newtheorem{Satz}{Theorem}[section]
\newtheorem{Prop}[Satz]{Proposition}
\newtheorem{Lem}[Satz]{Lemma}

\newtheorem{Thm}[Satz]{Theorem}
\newtheorem{Cor}[Satz]{Corollary}
\theoremstyle{definition}
\newtheorem{Dfn}[Satz]{Definition}
\newtheorem{Bem}[Satz]{Remark}

\newcommand{\tr}{\operatorname{Tr}}

\newcommand{\hess}{\operatorname{Hess}}

\newcommand{\vol}{{\operatorname{vol}}}
\newcommand{\dv}{\text{ }dv}

\newcommand{\norm}[1]{\left \lVert #1 \right \rVert}
\allowdisplaybreaks[1]

\renewcommand{\epsilon}{\varepsilon}

\newcommand{\R}{\ensuremath{\mathbb{R}}}

\newcommand{\s}{\ensuremath{\mathbb{S}}}
\numberwithin{equation}{section}

\usepackage{color}

\newtheorem{innercustomgeneric}{\customgenericname}
\providecommand{\customgenericname}{}
\newcommand{\newcustomtheorem}[2]{%
  \newenvironment{#1}[1]
  {%
   \renewcommand\customgenericname{#2}%
   \renewcommand\theinnercustomgeneric{##1}%
   \innercustomgeneric
  }
  {\endinnercustomgeneric}
}

\newcustomtheorem{customthm}{Theorem}
\newcustomtheorem{customcor}{Corollary}

\title{Stable proper biharmonic maps in Euclidean spheres}
\author{Volker Branding}
\date{\today}
\address{University of Vienna, Faculty of Mathematics\\
Oskar-Morgenstern-Platz 1, 1090 Vienna, Austria\\}
\email{volker.branding@univie.ac.at}

\author{Anna Siffert}
\address{Universität M\"unster, Mathematisches Institut\\
Einsteinstr. 62\\
48149 M\" unster\\
Germany}
\email{asiffert@uni-muenster.de}

\thanks{The first author gratefully acknowledges the support of the Austrian Science Fund (FWF) through the project "Geometric Analysis of Biwave Maps" (DOI: 10.55776/P34853). 
}

\begin{document}
\begin{abstract}
We construct an explicit family of stable proper weak biharmonic maps from the unit ball $B^m$, $m\geq 5$,
to Euclidean spheres.
To the best of the authors knowledge this is the first example of a stable proper weak biharmonic map from at compact domain.
To achieve our result we first establish the second
variation formula of the bienergy for maps from the unit ball into a Euclidean sphere. Employing this result, we examine the stability of
the proper weak biharmonic maps $q:B^m\rightarrow\s^{m^{\ell}}$, $m,\ell\in\mathbb{N}$ with $\ell\leq m$, which we recently constructed in \cite{BS25} and thus deduce the existence of an explicit family of stable proper biharmonic maps to Euclidean spheres.
\end{abstract} 
\maketitle

\section{Introduction}
The study of energy functionals is a central direction of research within geometric analysis. Such energy functionals can include various geometric objects such as maps between manifolds or sections of specific vector bundles. Typical questions arising here are the classification of critical points or understanding the properties of a given critical point as for example questions concerning its stability.

\smallskip

Recall that the \emph{energy} of a smooth map $\phi\colon M\to N$ between 
two Riemannian manifolds $(M,g)$ and $(N,h)$, is given by 
$$
E(\phi):=\frac{1}{2}\int_M |d \phi|^2 \dv,
$$
where $dv:=dv_g$ denotes the volume form of $M$ with respect to $g$.
Critical points of this energy functional are called \textit{harmonic maps} and they are characterized by the vanishing of the tension field $\tau$, which is defined by 
$$
\tau(\phi):= \operatorname{Tr}_g \bar\nabla d\phi.
$$
Here, $\bar\nabla$ denotes the connection on the pull-back bundle $\phi^{*}TN$.
The harmonic map equation \(\tau(\phi)=0\) comprises a semilinear elliptic partial differential equation of second order for which many results could be achieved over the last sixty years. For more details on harmonic maps one can consult the book \cite{MR2044031}.
\smallskip
In recent years many mathematicians got attracted by the study of higher order energy functionals which extend the energy of a map.
The most prominent generalization to mention here is
the (intrinsic) bienergy functional, or simply the bienergy, which is defined by
\begin{align}
\label{eq:bienergy}
E_2(\phi):=\frac{1}{2}\int_M |\tau (\phi)|^2 \dv.    
\end{align}

Critical points of the bienergy functional $E_2$ are called \textit{intrinsic biharmonic maps}, or simply \textit{biharmonic maps}, and
are characterized by the vanishing of the  \textit{bitension field} $\tau_2$, which is given by
\begin{align}
\label{tau2}    
\tau_2(\phi):= \bar\Delta \tau (\phi) +\textrm{Tr}_g R^N(\tau(\phi), d\phi) d\phi.
\end{align}
Here, $\bar\Delta$ denotes the rough Lapacian acting on sections of the pullback bundle $\phi^\ast TN$.

Trivially, harmonic maps are examples of biharmonic maps. 
Henceforth we will restrict our considerations to non-harmonic biharmonic maps, so called \textit{proper} biharmonic maps.
In the case that \(M\) is closed and \(N\) has non-positive curvature 
the maximum principle implies that all biharmonic maps are harmonic \cite{MR2640582}. In contrast, there exist many examples of proper biharmonic maps
in the Euclidean sphere, see for example \cite{MR4410183}.

For more details on biharmonic maps in Riemannian geometry we refer to the recent book \cite{MR4265170}.

\smallskip

In the specific case of maps to the sphere the  
bienergy \eqref{eq:bienergy} acquires the form
\begin{align*}
  E_2(u)=\frac{1}{2}\int_M\big(|\Delta u|^2-|\nabla u|^4)\dv, 
\end{align*}
where \(u\colon M\to\s^n\subset\R^{n+1}\), and whose critical points are given by
\begin{align}
\label{eq:biharmonic-intro}
\Delta^2u+2\operatorname{div}\big(|\nabla u|^2\nabla u\big)
-\big(\langle\Delta^2u,u\rangle-2|\nabla u|^4\big)u=0.
\end{align}

\smallskip

Although being natural from the geometric point of view, the 
bienergy \eqref{eq:bienergy} is not coercive and for this reason it is a challenging task to prove the existence of critical points. In this regard we want to point out the recent result of Laurain and Lin \cite{MR4227593} who showed the existence of 
a biharmonic map from the unit ball in four dimensions to the sphere using the heat flow method.

\smallskip

Recall that the
equator map is given by
\begin{align}
\label{equator}
u^{(1)}_{\star}\colon B^m&\rightarrow\s^m\subset\mathbb{R}^m\times\mathbb{R},\\ 
  \notag x&\mapsto\big(u^{(1)}(x),0\big),
\end{align}
where $r:=r(x):=\norm{x}$ denotes the Euclidean distance to the origin and $u^{(1)}:B^m\rightarrow\s^{m-1}, x\mapsto\frac{x}{r}$ is the well-studied radial projection map.
The equator map $u^{(1)}_{\star}$ is a well-known example of a harmonic map, see for example \cite{MR705882}, and thus does not provide an example of a proper biharmonic map.

Recently, Nakauchi \cite{MR4593065} constructed a generalized radial projection map $u^{(\ell)}:B^m\rightarrow\s^{m^\ell-1}$, where $\ell\in\mathbb{N}$ is a parameter and $u^{(1)}$ coincides with the radial projection map \(u^{(1)}\) introduced above. This gives rise to a generalized equator map as follows
\begin{align}
\label{equator-generalized}
u^{(\ell)}_{\star}\colon B^m&\rightarrow\s^{m^{\ell}}\subset\mathbb{R}^{m^{\ell}}\times\mathbb{R},\\ 
  \notag x&\mapsto\big(u^{(\ell)}(x),0\big),
\end{align}
which we studied extensively in \cite{BS25-2}.
It turns out that, by employing a suitable deformation, the generalized equator map \eqref{equator-generalized} can be turned into 
a biharmonic map, compare 
\cite[Theorem 3.2]{BS25}, \cite{MR4076824}.
More precisely, the map
\(q:B^m\rightarrow \s^{m^{\ell}}\)
given by
\begin{align}
\label{dfn:biharmonic-q}
q:=\big(\sin\alpha\cdot u^{(\ell)},\cos\alpha\big), \qquad \alpha\in (0,\frac{\pi}{2})
\end{align}
is a proper 
biharmonic map if and only if the following equation is satisfied 
\begin{align}
\label{eq:alpha}
\sin^2\alpha=\frac{\ell(\ell+m-2)+2m-8}{2\ell(\ell+m-2)}.
\end{align} 
It turns out that for \(m\geq 5,\ell\geq 2\) equation \eqref{eq:alpha}
always admits a solution and leads to a weak biharmonic map.
Note that these biharmonic maps would be smooth if we would consider \(B^m\setminus\{0\}\) as domain.

\bigskip

In this manuscript we will study a specific property of 
biharmonic maps, namely their stability:
A biharmonic map $u:M\rightarrow \s^n$ is called \textit{stable} if the second variation of \eqref{tau2}
is non-negative, i.e. the inequality
\begin{align}
\label{ineq:stable}
    \frac{d^2}{dt^2}E_2(u_t)\lvert_{t=0}\geq 0
\end{align}
holds for all variations $u_t$ of $u$ with $u_t-u\in W_0^{2,2}(M,\mathbb{R}^{n+1})$. 
If inequality (\ref{ineq:stable}) is strict for all variations $u_t$ of $u$ with $u_t-u\in W_0^{2,2}(M,\mathbb{R}^{n+1})$, the biharmonic map $u:M\rightarrow \s^n$ is called \textit{strictly stable}. 
If the map $u$ is not stable (and thus not strictly stable), it is called \textit{unstable}.

\smallskip

The main achievement of this manuscript is the construction 
of proper stable 
biharmonic maps in spheres.
In order to set up our main result we proceed as follows:
 First we establish the second variation of the bienergy for maps from compact manifolds with boundary to spheres. Employing this identity we deduce a sufficient condition for the stability of sphere-valued 
 biharmonic maps. We then apply this result to examine the stability of
the 
biharmonic maps $q:B^m\rightarrow\s^{m^{\ell}}$, $m,\ell\in\mathbb{N}$ with $\ell\leq m$, which we have constructed in \cite{BS25}
and recalled above.

Combining all these steps we finally arrive at the following result:

\begin{Satz}
\label{thm:intrinsic}
The proper 
biharmonic map \(q\colon B^m\to\s^{m^\ell},m\geq 5\) 
given by \eqref{dfn:biharmonic-q}
is 
strictly stable if 
the inequalities
\begin{align}
\label{eq:condition-ml-stable}
m> 2(\sqrt{12\ell^2+30\ell+12}+3\ell+6),\qquad l\leq m    
\end{align}
are satisfied.
\end{Satz}

Consequently, for each pair of integers $m,\ell\in\mathbb{N}$ such that  (\ref{eq:condition-ml-stable}) is satisfied and (\ref{eq:alpha}) admits a solution, we get an explicit strictly stable 
proper biharmonic map \(q:B^m\rightarrow \s^{m^{\ell}}\). 

\begin{Bem}
Let us connect the condition \eqref{eq:condition-ml-stable} to the results that have already been established in the literature:
\begin{enumerate}
    \item For \(\ell=1\) it is known that we obtain a biharmonic map for \(m=5,6\) as these are the only dimensions in which \eqref{eq:alpha} admits a solution, see \cite{MR4076824}.
Moreover, in that reference it was shown by constructing an explicit equivariant variation that these biharmonic maps are unstable.
This is consistent with Theorem \ref{thm:intrinsic} as for \(\ell=1,m=5,6\) inequality 
\eqref{eq:condition-ml-stable} is not satisfied: We would need
$m\geq 33$ in order for \eqref{eq:condition-ml-stable} to be satisfied when \(\ell=1\).
\item For $\ell=2$ the condition \eqref{eq:alpha} can be solved if \(m\geq 5\), see \cite[Theorem 1.2]{MR4830603}.
Moreover, in \cite[Theorem 1.3]{MR4830603} it was shown that this biharmonic map is unstable if \(5\leq m\leq 12\), again by constructing an explicit equivariant variation.
This is still consistent with \eqref{eq:condition-ml-stable} which,
for \(\ell=2\), yields stability if $m\geq 46$.

\item Moreover, for $\ell=3$ we also have that \eqref{eq:alpha} 
admits a solution for \(m\geq 5\), see \cite[Theorem 1.4]{MR4830603}.
This map is again unstable when \(5\leq m\leq 18\) being in agreement with
\eqref{eq:condition-ml-stable} which requires \(m\geq 59\)
in order to obtain a stable biharmonic map.

\item In \cite[Theorem 1.16]{MR4216418}
the first example of a strictly stable, proper biharmonic map into
a sphere has been provided by Montaldo, Oniciuc and Ratto:
Let $\varphi\colon\mathbb{R} \to \s^2$ be given by
\begin{equation}\label{A-harmonic-examples}
 \gamma \mapsto \, \left ( \cos(A(\gamma)),\,\sin(A(\gamma)), \,0\right ) \in \s^2\hookrightarrow  \R^3 \,,
\end{equation}
where
\begin{equation*}
A(\gamma)= a \gamma^3+b \gamma^2+c \gamma +d \,,
\end{equation*}
and $a,b,c,d\in\mathbb{R}$ with $a^2+b^2>0$. This provides a proper biharmonic curve on \(\s^2\) with a non-compact domain.
Assume that either
\begin{equation*}
 a=0 \quad {\rm or} \,\, \left\{a\neq 0 \,\,{\rm and} \,\, b^2 -3ac \leq 0 \right \}\,.
\end{equation*}
Then $\varphi$ is strictly stable with respect to compactly supported variations.

\item In \cite{MR3357596} 
biharmonic, rotationally symmetric, conformal maps between four-dimensional models of constant sectional curvature have been constructed. 
In particular, in \cite[Theorem 5.3]{MR3357596} the authors showed that
an explicitly given rotationally symmetric proper biharmonic conformal diffeomorphism \(\varphi\colon B^4\to\mathbb{H}^4\) is equivariantly stable.

Here, equivariant stability refers to the fact that only variations which are invariant under a group action are allowed, see e.g. \cite{MR4477489} and the references therein for more details.

One often sees the
phenomenon that a critical point of an energy functional is unstable in general but is stable with respect to equivariant variations. This can be explained by the fact that only a special class of directions in the formula for the second
variation of the energy is considered.

\item Concerning the stability of \eqref{equator} we want to mention the seminal
result of Jäger and Kaul \cite{MR705882} who showed that the equator map is minimizing if $m\geq 7$ and unstable if $3\leq m<7$ when considered as a harmonic map.
Later, in \cite{MR2322746} Hong and Thompson studied the stability of (\ref{equator}) as an extrinsic biharmonic map and showed that the equator map
is minimizing if $m\geq 10$ and unstable if $5\leq m<10$.
Recently, Fardoun, Montaldo and Ratto extended this stability analysis to the case of the equator map as a critical point of the extrinsic $k$-energy for $k\geq 3$
in \cite{MR4436204}. 
Very recently, we extended all of the above results to the case of the generalized equator map \eqref{equator-generalized} in \cite{BS25-2}.

\item For proper biharmonic maps which are given explicitly it is often possible to determine the complete spectrum of the corresponding Jacobi operator, see for example \cite{MR4216418}. However, these calculations are usually very demanding as one needs to solve a spectral problem for a Jacobi operator of fourth order.

\item In Theorem 18 of \cite{MR2640582} it is shown that a biharmonic map from a closed Riemannian manifold to the sphere is unstable under the assumption 
that it has constant energy density and satisfies the so-called conservation law.
As the map \eqref{dfn:biharmonic-q} does not have constant energy density and as the unit ball \(B^m\)
has a boundary the aforementioned result cannot be applied in our case.

\item Finally, we want to mention that all of the above results only seem to hold for the case of a compact domain manifold with boundary and that one cannot expect corresponding statements on closed manifolds, see  \cite{b25cag} for more details.

\item Biharmonic curves (biharmonic maps with a one-dimensional domain) from an interval with boundary have been studied by Mou in \cite{MR1774479}. 
\end{enumerate}
\end{Bem}

\bigskip

\textbf{Notation and Conventions:}
Throughout this article we will employ the following sign conventions: 
For the Riemannian curvature tensor field we use 
$$
R(X,Y)Z=[\nabla_X,\nabla_Y]Z-\nabla_{[X,Y]}Z,
$$ 
where \(X,Y,Z\) are vector fields.

For the rough Laplacian on the pull-back bundle $\phi^{\ast} TN$ we employ the analysts sign convention, i.e.
$$
\bar\Delta = \tr(\bar\nabla\bar\nabla-\bar\nabla_\nabla).
$$
In particular, this implies that the Laplace operator has a negative spectrum.

By $dv$ we represent the volume element of an arbitrary Riemannian manifold $(M,g)$. In case that the manifold \((M,g)\) has a non-empty boundary \(\partial M\) we denote the volume element on the boundary by
\(\dv^{\partial M}\).
Moreover, \(\nu\) will always denote the outward unit normal of \(\partial M\).
In the specific case that $(M,g)$ is the Euclidean ball, the volume element is simply written as $dx$.

Throughout this manuscript we employ the summation convention, i.e. we tacitly sum over repeated indices.
\medskip

\textbf{Organization:}
Section\,\ref{sec:prelim} contains preliminaries
on the generalized radial projection map introduced by Nakauchi
in \cite{MR4593065}.
In Section\,\ref{sec:boundary} the first and the second variation of the
bienergy for maps from domain manifolds with non-empty boundary are derived.
We provide stable proper 
biharmonic maps in Section\,\ref{sec:intrinsic}, in particular Theorem\,\ref{thm:intrinsic} is proven there.

\medskip

\textbf{Acknowledgements}: The authors would like to thank Florian Hanisch and Cezar Oniciuc for many inspiring comments on the manuscript.

\section{Preliminaries}
\label{sec:prelim}
In this section we provide basic facts on the generalized radial projection map introduced by Nakauchi
in \cite{MR4593065}, which appear in the definition (\ref{dfn:biharmonic-q}) of the map $q$. Here, we closely follow the presentation from \cite[Section 2]{BS25-2} 
and also show that the map \(q\colon B^m\to\s^{m^\ell}\) defined in \eqref{dfn:biharmonic-q} is a weak biharmonic map.

\medskip
\subsection{Generalized radial projection}
Recall that in \cite{MR4371934,MR4593065} a \lq generalized radial projection\rq\, has been introduced by Nakauchi.
More precisely, in \cite[Main Theorem, p.1]{MR4593065} Nakauchi showed that for any $\ell,m\in\mathbb{N}$ with $\ell\leq m$ there exists a harmonic map
\begin{align}
\label{nak-maps}
u^{(\ell)}\colon\mathbb{R}^{m}\setminus\{0\}&\to\s^{m^{\ell}-1},\\
\notag x=(x_1,\dots,x_m)&\mapsto u^{(\ell)}(x)=(u^{(\ell)}_{i_1\dots i_{\ell}}(x))_{1\leq i_1,\dots i_{\ell}\leq m}.
\end{align}
We set \(y_i=\frac{x_i}{r}\) and have the recursive definition
\begin{align*}
u^{(1)}_{i_1}(x)=&y_{i_1},\\
u^{(\ell)}_{i_1\ldots i_\ell}(x)=&C_{\ell,m}\big(y_{i_1}u^{(\ell-1)}_{i_1\ldots i_{\ell-1}}(x)-\frac{1}{\ell+m-3}r\frac{\partial}{\partial x_{i_\ell}}u^{(\ell-1)}_{i_1\ldots i_{\ell-1}}(x)\big),
\end{align*}
where 
\begin{align*}
 C_{\ell,m}=\sqrt{\frac{\ell+m-3}{2\ell+m-4}}.  
\end{align*}

In particular, we can use \eqref{nak-maps} to define a generalized equator map as 
we mentioned in the introduction, see \eqref{equator-generalized}.

For $\ell=1$ the equator map (\ref{equator}) is recovered.
Further, note that we consider $u_\star^{(\ell)}$ as one-parameter family of maps and use the singular to refer to it.

\smallskip

We will now give some more background on the map defined in \eqref{nak-maps}. Let $\ell,m\in\mathbb{N}$ such that $\ell\leq m$.
Nakauchi \cite{MR4593065} proved that the map $u^{(\ell)}\colon\mathbb{R}^{m}\setminus\{0\}\to\s^{m^\ell-1}$ has the following properties: 
\begin{enumerate}
\item $u^{(\ell)}$ satisfies the equation for harmonic maps to spheres
\begin{align*}
    \Delta u^{(\ell)}+\lvert\nabla u^{(\ell)}\rvert^2 u^{(\ell)}=0; 
\end{align*}
\item $u^{(\ell)}$ is a polynomial in $u_{i_1},\dots, u_{i_{\ell}}$ of degree $\ell$, where $u_{i_j}=\frac{x_{i_j}}{r}$;
\item $\lvert\nabla u^{(\ell)}\rvert^2=\frac{\ell(\ell+m-2)}{r^2}$.
\end{enumerate}
Furthermore, Nakauchi \cite[Proposition 1, (1)]{MR4593065} showed that the map $u^{(\ell)}$ satisfies the identity
\begin{align}
\label{ortho}
   \sum_{j=1}^m y_j\cdot\nabla_ju^{(\ell)}=0.
\end{align}

The following Lemma, which is \cite[Lemma 2.3]{BS25}, gives an explicit expression for $\Delta ^{k}u^{(\ell)}$ for any $k\in\mathbb{N}$.

\begin{Lem}
\label{lem:delta-nak}
 For each $k\in\mathbb{N}$ the map \(u^{(\ell)}\) defined in (\ref{nak-maps}) satisfies 
 \begin{align*}
          \Delta ^{k}u^{(\ell)}=
          \prod_{j=1}^k(2j+\ell-2)(2j-\ell-m)\frac{u^{(\ell)}}{r^{2k}}.
    \end{align*} 
\end{Lem}

\begin{Bem}
The equator map, i.e. $u_\star^{(1)}$, has been studied thoroughly in the last decades,
see e.g. \cite{MR4436204, MR705882}. 
Recently, the map $u^{(\ell)}$ with $\ell\geq 1$ has been a valuable tool to construct extrinsic polyharmonic maps and to study their stability: 
The map $u^{(\ell)}$ has been used in \cite{MR4830603,BS25} to provide proper biharmonic maps to spheres.
Further, the map $u^{(\ell)}$ has been used in \cite{BS25} to construct proper triharmonic maps to spheres.
A detailed stability analysis of the generalized equator map (\ref{equator-generalized}), considered as a critical point of both the extrinsic $k$-energy and the \(p\)-energy, has been provided in \cite{BS25-2}.
For the notions of extrinsic polyharmonic maps and the extrinsic $k$-energy we refer the reader to \cite{BS25,BS25-2} and the references therein.
\end{Bem}

\subsection{Weak biharmonic maps to spheres}
In many nonlinear problems in geometric analysis it is favorable to not only consider the class of smooth maps but to also allow for maps of lower regularity and to consider a weak version of the problem at hand.
In order to approach the notion of weak solutions of the harmonic map equation
let us recall the definition of the Sobolev space for maps to the sphere
\begin{align*}
W^{p,q}(M,\s^n):=\{u\in W^{p,q}(M,\R^{n+1})~~\mid u(x)\in\s^n ~~ a.e.\}.    
\end{align*}
In the case that \(\partial M\neq 0\) we let \(u_0\in W^{p,q}(M,\s^n)\) and define
\begin{align*}
W^{p,q}_{u_0}(M,\s^n):=\{u\in W^{p,q}(M,\s^n) \mid \nabla^k(u-u_0)\big|_{\partial M}=0,
0\leq k\leq p-1\}.
\end{align*}
Here, the boundary condition is to be understood in the sense of traces.

Now, we recall the following definition as it is central within this manuscript.
\begin{Dfn}
A map \(u\in W^{2,2}(B^m,\s^n)\) is called a weak biharmonic map if it solves \eqref{eq:biharmonic-intro} in the sense of distributions.
\end{Dfn}

In order to show that \(q\colon B^m\to \s^{m^\ell}\) is a weak biharmonic map we have to check if \(q\in W^{2,2}(B^m,\s^{m^\ell})\).
We find that
\begin{align*}
\int_{B^m}|\nabla q|^2\dv    
=&\frac{1}{2}\big(\ell(\ell+m-2)+2m-8\big)\vol(\s^{m-1})\int_0^1r^{m-3}dr,\\
\int_{B^m}|\Delta q|^2\dv    
=&\frac{1}{2}\ell(\ell+m-2)(\ell(\ell+m-2)+2m-8)\vol(\s^{m-1})\int_0^1r^{m-5}dr. 
\end{align*}
We realize that we need to require \(m\geq 5\) in order for the second integral to be finite,
hence \(q\) belongs to \(W^{2,2}(B^m,\s^{m^2-1})\) whenever \(m\geq 5\).

\section{Variational formulas for manifolds with boundary}
\label{sec:boundary}
In this section we derive the formulas for the first and second variation
of the 
bienergy on compact domain manifolds $M$ with non-empty, smooth boundary $\partial M$.

\smallskip

Throughout this section let \(\phi_t\colon (-\epsilon,\epsilon)\times M\to N,\epsilon>0\) be a smooth variation of the map \(\phi\colon M\to N\) with associated variational vector field
\begin{align}
\label{eq:variation-phi}
\frac{\bar\nabla\phi_t}{\partial t}\big|_{t=0}=V.
\end{align}

Moreover, \(\{e_i\},i=1,\ldots,\dim M\) is a local orthonormal frame of \(TM\)
which, at a fixed but arbitrary point \(p\in M\) satisfies
\(\nabla_{e_j}e_i=0,1\leq i,j\leq\dim M\) and
\(\nabla_{\partial_t}e_j=0,1\leq j\leq\dim M\).

Throughout the manuscript we will frequently make use of the following Green's formula
for the rough Laplacian. 
\begin{Lem}
Let \((M,g)\) be a compact Riemannian manifold with non-empty boundary \(\partial M\) and \(E\) a vector bundle over \(M\). Then, the following
formula holds
\begin{align}
\label{eq:green}
\int_M\langle\Delta^E V,W\rangle\dv
-\int_M\langle V,\Delta^E W\rangle\dv
=\int_{\partial M}\big(
\langle\nabla^E_{\nu}V,W\rangle-\langle V,\nabla^E_\nu W\rangle
\big)\dv^{\partial M},
\end{align}
where \(V,W\in\Gamma(E)\), \(\nu\) represents the unit outward normal of \(\partial M\)
and \(\Delta^E,\nabla^E\) the rough Laplacian and the connection on \(E\) respectively.
\end{Lem}

\begin{proof}
This follows from 
\begin{align*}
\langle V,\Delta^EW\rangle+\langle\nabla^EV,\nabla^EW\rangle=
\operatorname{div}\langle V,\nabla^E W\rangle
\end{align*}
and the divergence theorem
\begin{align*}
\int_M\operatorname{div}\langle V,\nabla^E W\rangle\dv
=\int_{\partial M}\langle V,\nabla_{\nu}^E W\rangle\dv^{\partial M},
\end{align*}
see \cite[Section 2.2]{MR2044031} for more details.
\end{proof}

\begin{Thm}[First variation formula with boundary]
\label{prop-bound}
Consider a variation as described above. Then, the first variation of
the bienergy functional acquires the form
\begin{align}
\label{eq:first-boundary}
\frac{d}{dt}\frac{1}{2}\int_M|\tau(\phi_t)|^2\dv
=&\int_M\langle\tau_2(\phi_t),d\phi_t(\partial_t)\rangle\dv \\
\nonumber&+\int_{\partial M}\big(
\langle\bar\nabla_\nu d\phi_t(\partial_t),\tau(\phi_t)\rangle-\langle d\phi_t(\partial_t),\bar\nabla_\nu\tau(\phi_t)\rangle
\big)\dv^{\partial M}.
\end{align}
\end{Thm}

\begin{proof}
We calculate
\begin{align*}
\frac{d}{dt}\frac{1}{2}\int_M|\tau(\phi_t)|^2\dv=&
\int_M\langle R^N(d\phi_t(\partial_t),d\phi_t(e_j))d\phi_t(e_j),\tau(\phi_t)
\rangle\dv
+\int_M\langle\bar\Delta d\phi_t(\partial_t),\tau(\phi_t)\rangle\dv.
\end{align*}
Now, we apply Green's formula for manifolds with boundary, that is 
\eqref{eq:green} and find
\begin{align*}
\int_M\langle\bar\Delta d\phi_t(\partial_t),\tau(\phi_t)\rangle\dv
=&\int_M\langle d\phi_t(\partial_t),\bar\Delta\tau(\phi_t)\rangle\dv \\
&+\int_{\partial M}\big(\langle\bar\nabla_\nu d\phi_t(\partial_t),\tau(\phi_t)\rangle-\langle d\phi_t(\partial_t),\bar\nabla_\nu\tau(\phi_t)\rangle
\big)\dv^{\partial M}.
\end{align*}
Combining these two identities and using (\ref{tau2}) completes the proof.
\end{proof}

A direct consequence of Theorem\,\ref{prop-bound} is the following fact:
\begin{Prop}
Let \((M,g)\) be a compact Riemannian manifold with non-empty boundary \(\partial M\) and \((N,h)\) a Riemannian manifold. Then, a smooth map \(\phi\colon M\to N\) is biharmonic if it satisfies
\begin{align}
\label{eq:bc}
\begin{cases}
    \tau_2(\phi)=0 \quad&\text{ on } M, \\
    \phi =\xi \quad&\text{ on }\partial M,\\
    d\phi(\partial_\nu) = \chi\quad&\text{ on } \partial M,
\end{cases}
\end{align}
where \(\xi\colon\partial M\to N,\chi\colon M\to TN\) are given smooth boundary data.
\end{Prop}

\begin{Thm}
\label{prop:second}
Let \((M,g)\) be a compact Riemannian manifold with non-empty boundary
\(\partial M\neq \emptyset\) and \((N,h)\) a second Riemannian manifold.
Assume that \(\phi\colon M\to N\) is a smooth biharmonic map satisfying
the boundary conditions \eqref{eq:bc}. Then, the second variation of the bienergy acquires the form
\begin{align}
\label{eq:sv-intrinsic}
\hess E_2(\phi)(V,V)=&\int_M|\bar\Delta V-R^N(d\phi(e_i),V)d\phi(e_i)|^2\dv \\
\nonumber&-\int_M\langle (\nabla_{d\phi(e_i)}R^N)(d\phi(e_i),\tau(\phi))V,V\rangle\dv \\
\nonumber&-\int_M\langle (\nabla_{\tau(\phi)}R^N)(d\phi(e_i),V)d\phi(e_i),V\rangle\dv \\
\nonumber&-\int_M\langle R^N(\tau(\phi),V)\tau(\phi),V\rangle\dv\\
\nonumber&-2\int_M\langle R^N(d\phi(e_i),V)\bar\nabla_{e_i}\tau(\phi),V\rangle\dv\\
\nonumber&-2\int_M\langle R^N(d\phi(e_i),\tau(\phi))\bar\nabla_{e_i}V,V\rangle\dv,
\end{align}
where \(V\in\Gamma(\phi^\ast TN)\) represents the variational vector field.
\end{Thm}

\begin{proof}
We consider a variation of the map \(\phi\colon M\to N\) as detailed in \eqref{eq:variation-phi}. Then, we find
\begin{align*}
\frac{\bar\nabla}{\partial t}\big|_{t=0}\bar\Delta\tau(\phi)
=&R^N(V,d\phi(e_i))\bar\nabla_{e_i}\tau(\phi)
+\bar\nabla_{e_i}\big(R^N(V,d\phi(e_i))\tau(\phi)\big) \\
&+\bar\Delta\big(\bar\Delta V-R^N(d\phi(e_i),V)d\phi(e_i)\big) \\
=&2R^N(V,d\phi(e_i))\bar\nabla_{e_i}\tau(\phi)
+(\nabla_{d\phi(e_i)}R^N)(V,d\phi(e_i))\tau(\phi) \\
&+R^N(\bar\nabla_{e_i}V,d\phi(e_i))\tau(\phi)
+R^N(V,\tau(\phi))\tau(\phi)\\
&+\bar\Delta\big(\bar\Delta V-R^N(d\phi(e_i),V)d\phi(e_i)\big),
\end{align*}
where we used that
\begin{align}
\label{eq:id-t-tension}
\frac{\bar\nabla}{\partial t}\big|_{t=0}\tau(\phi)
=\bar\Delta V-R^N(d\phi(e_i),V)d\phi(e_i)
\end{align}
in the first step. Moreover, we get
\begin{align*}
 \frac{\bar\nabla}{\partial t}\big|_{t=0} R^N(\tau(\phi),d\phi(e_i))d\phi(e_i)   
 =&(\nabla_VR^N)(\tau(\phi),d\phi(e_i))d\phi(e_i)
 +R^N(\bar\Delta V,d\phi(e_i))d\phi(e_i) \\
 &-R^N(R^N(d\phi(e_k),V)d\phi(e_k),d\phi(e_i))d\phi(e_i) \\
 &+R^N(\tau(\phi),\bar\nabla_{e_i}V)d\phi(e_i)
 +R^N(\tau(\phi),d\phi(e_i))\bar\nabla_{e_i}V \\
=&-(\nabla_{d\phi(e_i)}R^N)(V,\tau(\phi))d\phi(e_i)
-(\nabla_{\tau(\phi)}R^N)(d\phi(e_i),V)d\phi(e_i)
\\
&+R^N(\bar\Delta V,d\phi(e_i))d\phi(e_i) \\
 &-R^N(R^N(d\phi(e_k),V)d\phi(e_k),d\phi(e_i))d\phi(e_i) \\
 &+R^N(\tau(\phi),\bar\nabla_{e_i}V)d\phi(e_i)
 +R^N(\tau(\phi),d\phi(e_i))\bar\nabla_{e_i}V,
\end{align*}
where we used \eqref{eq:id-t-tension} to establish the first equality.
For the second equality we employed the first and the second Bianchi identity.

To determine the boundary contributions in the second variation 
formula of the bienergy we calculate
\begin{align*}
\frac{\bar\nabla}{\partial t}\big|_{t=0}    
\langle d\phi_t(\partial_t),\bar\nabla_\nu\tau(\phi_t)\rangle=&
\langle \frac{\bar\nabla}{\partial t}d\phi_t(\partial_t)\big|_{t=0},
\bar\nabla_\nu\tau(\phi)\rangle 
+\langle R^N(V,d\phi(\partial_\nu))\tau(\phi),V\rangle \\
&+\langle V,\bar\nabla_\nu\bar\Delta V\rangle
-\langle V,\bar\nabla_\nu\big(R^N(d\phi(e_i),V)d\phi(e_i)\big)\rangle
\end{align*}
and 
\begin{align*}
\frac{\bar\nabla}{\partial t}\big|_{t=0}    
\langle\bar\nabla_\nu d\phi_t(\partial_t),\tau(\phi)\rangle
=&\langle \frac{\bar\nabla}{\partial t}\bar\nabla_\nu d\phi_t(\partial_t)\big|_{t=0},\tau(\phi_t)\rangle \\
&+\langle\bar\nabla_\nu V,\bar\Delta V\rangle
-\langle\bar\nabla_\nu V,R^N(d\phi(e_i),V)d\phi(e_i)\rangle.
\end{align*}

Thus, we arrive at
\begin{align*}
\hess E_2(\phi)(V,V)=&\int_M\langle \frac{\bar\nabla}{\partial t}d\phi_t(\partial_t)\big|_{t=0},\tau_2(\phi)\rangle\dv \\
&-\int_{\partial M}\langle \frac{\bar\nabla}{\partial t}d\phi_t(\partial_t)\big|_{t=0},
\bar\nabla_\nu\tau(\phi)\rangle\dv^{\partial M} 
+\int_{\partial M}\langle \frac{\bar\nabla}{\partial t}\bar\nabla_\nu d\phi_t(\partial_t)\big|_{t=0},\tau(\phi)\rangle\dv^{\partial M}
\\
&+\int_M\langle V,\bar\Delta(\bar\Delta V-R^N(d\phi(e_i),V)d\phi(e_i))\rangle\dv \\
&+\int_M\langle R^N(\bar\Delta V,d\phi(e_i))d\phi(e_i),V\rangle\dv\\
&+\int_M|R^N(d\phi(e_i),V)d\phi(e_i)|^2\dv\\
&+\int_M\langle R^N(V,\tau(\phi))\tau(\phi),V\rangle\dv \\
&+2\int_M\langle R^N(V,d\phi(e_i))\bar\nabla_{e_i}\tau(\phi),V\rangle\dv \\
&+2\int_M\langle R^N(\tau(\phi),d\phi(e_i))\bar\nabla_{e_i}V,V\rangle\dv \\
&+\int_M\langle(\nabla_{d\phi(e_i)}R^N)(\tau(\phi),d\phi(e_i))V,V\rangle\dv\\
&-\int_M\langle(\nabla_{\tau(\phi)}R^N)(d\phi(e_i),V)d\phi(e_i),V\rangle\dv \\
&-\int_{\partial M}\langle R^N(V,d\phi(\partial_\nu))\tau(\phi),V\rangle \dv^{\partial M}\\
&-\int_{\partial M}\langle\bar\nabla_\nu\bar\Delta V,V\rangle\dv^{\partial M}
+\int_{\partial M}\langle\bar\nabla_\nu\big(R^N(d\phi(e_i),V)d\phi(e_i)\big),V\rangle\dv^{\partial M} \\
&+\int_{\partial M}\langle\bar\nabla_\nu V,\bar\Delta V\rangle\dv^{\partial M}
-\int_{\partial M}\langle\bar\nabla_\nu V,R^N(d\phi(e_i),V)d\phi(e_i)\rangle
\dv^{\partial M},
\end{align*}
where we used the first Bianchi identity as follows
\begin{align*}
&R^N(\tau(\phi),\bar\nabla_{e_i}V)d\phi(e_i)
+R^N(\tau(\phi),d\phi(e_i))\bar\nabla_{e_i}V
+R^N(\bar\nabla_{e_i}V,d\phi(e_i))\tau(\phi) \\
=&2R^N(\tau(\phi),d\phi(e_i))\bar\nabla_{e_i}V, \\
&-(\nabla_{d\phi(e_i)}R^N)(V,\tau(\phi))d\phi(e_i)
+(\nabla_{d\phi(e_i)}R^N)(V,d\phi(e_i))\tau(\phi) 
=(\nabla_{d\phi(e_i)}R^N)(\tau(\phi),d\phi(e_i))V.
\end{align*}
Concerning the fourth term in the above equation for $\hess E_2(\phi)(V,V)$ we again use 
Green's formula (\ref{eq:green}) to deduce
\begin{align*}
\int_M\langle V,\bar\Delta(\bar\Delta V-R^N(d\phi(e_i),V)d\phi(e_i))\rangle\dv
=&\int_M\langle \bar\Delta V,\bar\Delta V-R^N(d\phi(e_i),V)d\phi(e_i)\rangle\dv\\
&-\int_{\partial M}\langle\bar\nabla_\nu V,\bar\Delta V-R^N(d\phi(e_i),V)d\phi(e_i))\rangle \dv^{\partial M} \\
&+\int_{\partial M}\langle V,\bar\nabla_\nu(\bar\Delta V-R^N(d\phi(e_i),V)d\phi(e_i))\rangle \dv^{\partial M}.
\end{align*}
The claim now follows from adding up all the contributions and using
that \(\phi\colon M\to N\) is a biharmonic map satisfying the boundary data.
\end{proof}

\begin{Bem}
If $\partial M=\emptyset$, the formula for the second variation \eqref{eq:sv-intrinsic} agrees with the formula provided by Jiang \cite[Theorem 16]{MR2640582} who studied the case of a closed domain manifold only. Note that a different sign convention for the curvature tensor is used in \cite{MR2640582}.
\end{Bem}

\section{Stable biharmonic maps}
\label{sec:intrinsic}
In this section we employ the generalized radial projection \eqref{nak-maps} to obtain the first example of a strictly stable, proper 
biharmonic map, i.e. we prove Theorem \ref{thm:intrinsic}. 

We proceed as follows: First, we establish the second variation of the bienergy for sphere-valued maps. From this identity we deduce a sufficient condition for the stability of sphere-valued  
biharmonic maps. 
In a second step we then apply this result to examine the stability of
the sphere-valued 
biharmonic maps $q:B^m\rightarrow\s^{m^{\ell}}$, $m,\ell\in\mathbb{N}$, which we recently constructed in \cite{BS25}. 

Combining these results we can then deduce
the existence of strictly stable, proper 
biharmonic maps
which are given explicitly.

\smallskip

As already mentioned in the introduction, in the specific case of maps to the sphere \(u\colon M\to\s^n\), the bienergy $E_2$ acquires the form
\begin{align}
\label{eq:bienergy-u-sphere}
  E_2(u)=\frac{1}{2}\int_M\big(|\Delta u|^2-|d u|^4)\dv.  
\end{align}
First, we recall the precise structure of the critical points of \eqref{eq:bienergy-u-sphere}.

\begin{Prop}
Let \((M,g)\) be a compact Riemannian manifold with non-empty boundary \(\partial M\). Then, a smooth map \(u\colon M\to\s^n\subset\R^{n+1}\) is biharmonic if it satisfies
\begin{align}
\label{eq:bc-u}
\begin{cases}
    \tau_2(u)=0 \quad&\text{ on } M, \\
    u =\xi \quad&\text{ on }\partial M, \\
    du(\partial_\nu) = \chi\quad&\text{ on } \partial M,
\end{cases}
\end{align}
where \(\xi\colon\partial M\to\s^n,\chi\colon\partial M\to T\s^n\) are given smooth boundary data and \(\tau_2(u)\) is defined by
\begin{align*}
\tau_2(u)=\Delta^2u+2\operatorname{div}\big(|du|^2du\big)
-\big(\langle\Delta^2u,u\rangle-2|du|^4\big)u.   
\end{align*}
\end{Prop}

In the following Lemma we provide the second variation of the functional (\ref{eq:bienergy-u-sphere}) on Euclidean balls:

\begin{Lem}
\label{lem:second-var}
Let
\(u\colon B^m\to\s^n\subset\R^{n+1}\) be a solution of \eqref{eq:bc-u}.
Then, the Hessian of the bienergy is given by
\begin{align}
\label{eq:hess-bienergy-sphere}
\frac{d^2}{dt^2}\big|_{t=0}\frac{1}{2}\int_{B^m}\big(|\Delta u_t|^2
-|d u_t|^4\big)\,dx=&\int_{B^m}\big(|\Delta\eta|^2+2|du|^4|\eta|^2\big)\,dx    -\int_{B^m}\langle\Delta^2u,u\rangle|\eta|^2\,dx \\
\nonumber &
-2\int_{B^m}|du|^2|d\eta|^2\, dx-4\int_{B^m}|\langle du,d\eta\rangle|^2\, dx,
\end{align}
where \(\eta=\frac{du_t}{dt}\big|_{t=0}\in C^\infty(B^m,\R^{n+1})\) with \(\langle \eta,u\rangle=0\).
\end{Lem}

\begin{proof}
This follows from \cite[Proof of Proposition 2.2]{MR4436204} and the application of the second variation formula for the \(p\)-energy \eqref{eq:sv-p-energy} with \(p=4\): 
In order to determine the second variation of $E_2$, we calculate 
the second variations of the functionals $\frac{1}{2}\int_{B^m}|\Delta u|^2 \,dx$ and $-\frac{1}{2}\int_{B^m}|du|^4\,dx$ separately.

\smallskip

The second variation of the first term of \eqref{eq:bienergy-u-sphere} was determined 
in the proof of Proposition 2.2 in \cite{MR4436204}, namely it was shown  that
\begin{align}
\label{func-1}
\frac{d^2}{dt^2}\big|_{t=0}\frac{1}{2}\int_{B^m}|\Delta u_t|^2
\,dx=&\int_{B^m}|\Delta\eta|^2\,dx -\int_{B^m}\langle\Delta^2u,u\rangle|\eta|^2\,dx.
\end{align}

The second term in \eqref{eq:bienergy-u-sphere} corresponds, up to a constant factor, to the $p$-energy with $p=4$.
Therefore we recall that the second variation formula of the \(p\)-energy
\begin{align*}
E_p(u):=\frac{1}{p}\int_{B^m}|du|^p\, dx
\end{align*}
for maps from the ball to the sphere, which is given by
\begin{align}
\label{eq:sv-p-energy}
\frac{d^2}{dt^2}E_p(u_t)\big|_{t=0}=\int_{B^m}|du|^{p-2}\big(|d\eta|^2-|du|^2|\eta|^2\big)\, dx+(p-2)\int_{B^m}|du|^{p-4}|\langle du,d\eta\rangle|^2 \,dx,
\end{align}
 where \(\eta=\frac{du_t}{dt}\big|_{t=0}\) and 
\(\langle\eta,u\rangle=0\) as \(\eta\in\Gamma(T\s^n)\).
A derivation of 
\eqref{eq:sv-p-energy} can be found in \cite[p. 143]{MR1619840}. 

\smallskip

The claim follows immediately from identities (\ref{func-1}) and (\ref{eq:sv-p-energy}).
\end{proof}

Before we proceed with the construction of a proper stable biharmonic map, we will show in the next Lemma that, for the specific setting under consideration,
the formula for the second variation of the bienergy \eqref{eq:hess-bienergy-sphere} coincides with the one derived
by Jiang \cite{MR2640582}, which we recalled in Theorem \ref{prop:second}.
Let \(\phi\colon M\to N\) be a biharmonic map between Riemannian manifolds and
assume that \(N\) has constant curvature. In this case the second variation of the bienergy evaluated at a critical point \eqref{eq:sv-intrinsic} acquires the form
\begin{align}
\label{eq:sv-bienergy-jiang}
\hess E_2(\phi)(V,V)=&\int_M\big(
|\bar\Delta V-R^N(d\phi(e_i),V)d\phi(e_i)|^2
-\langle R^N(\tau(\phi),V)\tau(\phi),V\rangle\\
\nonumber&-2\langle R^N(d\phi(e_i),V)\bar\nabla_{e_i}\tau(\phi),V\rangle
-2\langle R^N(d\phi(e_i),\tau(\phi))\bar\nabla_{e_i}V,V\rangle
\big)
\dv ,
\end{align}
where \(V\in\Gamma(\phi^\ast TN)\).
 Note that in 
\cite{MR2640582} a different sign convention for the Riemann curvature tensor is used.

\begin{Lem}
Let 
\(\phi\colon B^m\to \s^n\) be a biharmonic map satisfying the boundary conditions, i.e. a solution of \eqref{eq:bc}.
Then
equation \eqref{eq:sv-bienergy-jiang}
is equivalent to equation \eqref{eq:hess-bienergy-sphere}.
\end{Lem}

\begin{proof}
Since \(N=\s^n\), the Riemann curvature tensor $R^N$ is given by  
\begin{align*}
R^N(X,Y)Z=\langle Y,Z\rangle X-\langle X,Z\rangle Y.
\end{align*}
Thus, \eqref{eq:sv-bienergy-jiang} acquires the form 
\begin{align*}
\hess E_2(\phi)(V,V)=&\int_{B^m}\big(
\big|\bar\Delta V+|d\phi|^2V-\langle V,d\phi(e_j)\rangle d\phi(e_j)|^2 \\
&-|\langle V,\tau(\phi)\rangle|^2+|\tau(\phi)|^2|V|^2
-2\langle V,\bar\nabla_{e_i}\tau(\phi)\rangle\langle d\phi(e_i),V\rangle \\
&+2\langle d\phi(e_i),\bar\nabla_{e_i}\tau(\phi)\rangle|V|^2
-2\langle\tau(\phi),\bar\nabla_{e_i}V\rangle\langle d\phi(e_i),V\rangle \\
&+2\langle d\phi(e_i),\bar\nabla_{e_i}V\rangle\langle\tau(\phi),V\rangle
\big)\,dx.
\end{align*}
Now, for \(u:=\iota\circ\phi\colon B^m\to\s^n\subset\R^{n+1}\),
with \(\iota\) being the canonical embedding of \(\s^n\) into \(\R^{n+1}\),
we have \(\tau(u)=\Delta u+|du|^2u\) and
\begin{align*}
\bar\nabla_{e_j}V&=\nabla_{e_j}V+\langle du(e_j),V\rangle u,\\  
\bar\Delta V&=\Delta V-\langle u,\Delta V\rangle u
+\langle du, V\rangle du,
\end{align*}
where we used that \(\langle u,V\rangle =0\). In the following we will not distinguish between \(V\) and \(d\iota(V)\) in order to shorten the presentation.

Furthermore, direct computations yield the following identities
\begin{align*}
|\langle V,\tau(\phi)\rangle|^2&=|\langle V,\Delta u\rangle|^2
,\\
|\tau(\phi)|^2|V|^2&=(|\Delta u|^2-|du|^4)|V|^2
,\\
\langle V,\bar\nabla_{e_j}\tau(\phi)\rangle\langle d\phi(e_j),V\rangle&=
(\langle V,\nabla\Delta u\rangle+\langle V, du\rangle |du|^2)\langle du,V\rangle
,\\
\langle d\phi(e_j),\bar\nabla_{e_j}\tau(\phi)\rangle |V|^2&=    
\langle\nabla u,\nabla\Delta u\rangle |V|^2+|du|^4|V|^2
,\\
\langle\tau(\phi),\bar\nabla_{e_j}V\rangle\langle d\phi(e_j),V\rangle&=
(\langle\Delta u,\nabla V\rangle-\langle du,V\rangle|du|^2)\langle du,V\rangle
,\\
\langle d\phi(e_j),\bar\nabla_{e_j}V\rangle\langle V,\tau(\phi)\rangle&=
\langle V,\Delta u\rangle\langle du,\nabla V\rangle
\end{align*}
such that
\begin{align*}\hess E_2(u)(V,V)=\int_{B^m}\big(&
|\Delta V|^2+2|du|^4|V|^2+(|\Delta u|^2+2\langle du,d\Delta u\rangle)|V|^2\\
&-|\langle u,\Delta V\rangle|^2 -|\langle V,\Delta u\rangle|^2 
+2|du|^2\langle\Delta V,V\rangle \\
&-2\langle V,\nabla\Delta u\rangle\langle du,V\rangle
-2\langle\Delta u,\nabla V\rangle\langle du,V\rangle
+2\langle V,\Delta u\rangle\langle du,\nabla V\rangle
\big)\,dx.
\end{align*}
Now, we employ the following equations
\begin{align*}
0=&\langle\Delta u,V\rangle +2\langle du,\nabla V\rangle
+\langle u,\Delta V\rangle, \\
0=&\Delta|du|^2+|\Delta u|^2+2\langle d\Delta u,du\rangle
+\langle\Delta^2u,u\rangle
\end{align*}
which follow from differentiating \(\langle u,V\rangle=0\) and 
\(|u|^2=1\).
Hence, we get
\begin{align*}
\hess E_2(u)(V,V)=\int_{B^m}\big(&
|\Delta V|^2+2|du|^4|V|^2-\langle\Delta^2u,u\rangle|V|^2
-4|\langle du,\nabla V\rangle|^2\\
&-|\langle u,\Delta V\rangle|^2 -|\langle V,\Delta u\rangle|^2 
-(\Delta|du|^2)|V|^2+2|du|^2\langle\Delta V,V\rangle \\
&-2\langle V,\nabla\Delta u\rangle\langle du,V\rangle
-2\langle\Delta u,\nabla V\rangle\langle du,V\rangle
-2\langle \Delta V,u\rangle\langle du,\nabla V\rangle
\big)\,dx.
\end{align*}
Using Green's formula for the Laplacian we obtain
\begin{align*}
\int_{B^m}\big(2|du|^2\langle\Delta V,V\rangle-(\Delta|du|^2)|V|^2\big)\dv
=&-2\int_{B^m}|du|^2|\nabla V|^2\, dx 
-\int_{\partial B^m}|V|^2\partial_{\nu}|du|^2\dv^{\partial B^m} \\
&+\int_{\partial B^m}|du|^2\partial_{\nu}|V|^2\dv^{\partial B^m}.
\end{align*}
Thus, it remains to show that
\begin{align}
\label{eq:identity-a}
0=\int_{B^m}\big(&
-|\langle u,\Delta V\rangle|^2 -|\langle V,\Delta u\rangle|^2 
\\
\nonumber&-2\langle V,d\Delta u\rangle\langle du,V\rangle
-2\langle\Delta u,\nabla V\rangle\langle du,V\rangle
-2\langle \Delta V,u\rangle\langle du,\nabla V\rangle
\big)\,dx.
\end{align}
Again using the divergence theorem we get
\begin{align*}
\int_{B^m}\big(\langle V,\nabla\Delta u\rangle\langle du,V\rangle
+\langle\Delta u,\nabla V\rangle\langle du,V\rangle\big)\,dx
=&\int_{B^m}\big(-|\langle V,\Delta u\rangle|^2
-\langle V,\Delta u\rangle\langle du,\nabla V\rangle
\big)\dv \\
&+\int_{\partial B^m}\langle V,\Delta u\rangle\langle du(\partial_\nu),V\rangle \dv^{\partial B^m}.
\end{align*}
Moreover, we have
\begin{align*}
-|\langle u,\Delta V\rangle|^2 -|\langle V,\Delta u\rangle|^2     
=-2|\langle V,\Delta u\rangle|^2
-4\langle\nabla u,\nabla V\rangle\langle\Delta u,V\rangle    
-4|\langle\nabla u,\nabla V\rangle|^2
\end{align*}
such that \eqref{eq:identity-a}
indeed holds true completing the proof.
\end{proof}

We 
now turn to the construction of stable proper biharmonic maps to the Euclidean sphere.

More precisely, from Lemma\,\ref{lem:second-var} we deduce the following sufficient condition for the stability of a sphere-valued biharmonic map:

\begin{Cor}
\label{cor:stability}
A
biharmonic map \(u\colon B^m\to\s^n\subset\R^{n+1}\) 
satisfying the boundary conditions
is strictly stable if the inequality
\begin{align}
\label{eq:stability-intrinsic}
\int_{B^m}\big(|\Delta\eta|^2+2|du|^4|\eta|^2\big)\,dx    
-\int_{B^m}\langle\Delta^2u,u\rangle|\eta|^2\,dx 
-6\int_{B^m}|du|^2|d\eta|^2\,dx >0
\end{align}
holds, where \(\eta=\frac{du_t}{dt}\big|_{t=0}\)  with \(\langle \eta,u\rangle=0\).
\end{Cor}

\begin{proof}
The Cauchy-Schwarz inequality yields
\begin{align*}
-2\int_{B^m}|du|^2|d\eta|^2\,dx-4\int_{B^m}|\langle du,d\eta\rangle|^2\,dx
\geq &-6 \int_{B^m}|du|^2|d\eta|^2\,dx.
\end{align*}
Plugging this inequality into the stability condition \eqref{eq:hess-bienergy-sphere} provides the claim.
\end{proof}

We now return to the construction of
stable biharmonic maps. More precisely, we will
use Corollary\ref{cor:stability} to examine the stability of specific biharmonic maps to spheres which we recently constructed in \cite{BS25}.
To this end, we first recall this existence result, see \cite[Theorem 3.2]{BS25}:

\begin{Satz}
The map
\(q:B^m\to \s^{m^{\ell}}\)
given by
\begin{align}
\label{dfq:biharmonic-q}
q:=\big(\sin\alpha\cdot u^{(\ell)},\cos\alpha\big),
\end{align}
where \(u^{(\ell)}\) is defined in \eqref{nak-maps},
is a proper weak biharmonic map if and only if the following equation is satisfied 
\begin{align*}
\sin^2\alpha=\frac{\ell(\ell+m-2)+2m-8}{2\ell(\ell+m-2)},\qquad \ell\leq m.
\end{align*}     
\end{Satz}

Note that \cite[Theorem 3.2]{BS25} is formulated for maps from \(\R^m\setminus\{0\}\) but, apart from including the boundary conditions, 
all arguments work the same way when considering \(B^m\) instead of \(\R^m\setminus\{0\}\) and changing to the weak formulation of the biharmonic map equation.

Straightforward calculations yield the following identities for (\ref{dfq:biharmonic-q}):
\begin{align}
\label{prop-q}
|dq|^2=&\frac{1}{2}\big(\ell(\ell+m-2)+2m-8\big)\frac{1}{r^2},\\
\int_{B^m}|dq|^4|\eta|^2\,dx=\notag&\frac{1}{4}\big(\ell(\ell+m-2)+2m-8\big)^2\int_{B^m}\frac{|\eta|^2}{r^4}\,dx, \\
\int_{B^m}\langle\Delta^2q,q\rangle|\eta|^2\,dx=\notag&
\frac{1}{2}\big(\ell+2\big)\big(m+\ell-4\big)
\big(\ell(\ell+m-2)+2m-8\big)\int_{B^m}\frac{|\eta|^2}{r^4}\,dx,\\
\int_{B^m}|dq|^2|d\eta|^2\,dx=\notag&\frac{1}{2}\big(\ell(\ell+m-2)+2m-8\big)\int_{B^m}\frac{|d\eta|^2}{r^2}\,dx.
\end{align}

Our main tool will be the following Hardy-type inequality:
\begin{Lem}
For \(\eta\in W^{1,2}_0(B^m,\s^n)\) with \(m\geq 5\) the following Hardy-type inequality holds
\begin{align}
\label{eq:hardy-laplace}
\int_{B^m}|\Delta\eta|^2\,dx\geq 
\frac{m^2}{4}\int_{B^m}\frac{|d\eta|^2}{r^{2}}\,dx.
\end{align}
\end{Lem}

\begin{proof}
The statement follows from \cite[Theorem 1.1]{MR2296305}
and a summation over the components, see \cite[p. 272]{MR4817500} for a similar discussion.  
\end{proof}

\begin{Bem}
There exist various refined versions of the Hardy inequality \eqref{eq:hardy-laplace} in the literature but \eqref{eq:hardy-laplace}
is sufficient for our purposes.
\end{Bem}

In the next lemma we show that there do not exist parallel vector fields on the spheres \(\s^n,n\geq 2\). This result is well-known, see e.g.
the manuscripts \cite{MR192436, MR451257, MR1079056} from which one can deduce this statement. However, we could not find a reference in the literature in which this result is formulated.

\begin{Lem}
\label{lem:parallel_vf}
There do not exist parallel vector fields on \(\s^n,n\geq 2\).    
\end{Lem}

\begin{proof}
Let \(Y\) be a parallel vector field on \(\s^n\), i.e. \(\nabla Y=0\), and \(X\) a vector field
that is not a constant multiple of \(Y\). Then, we have
\begin{align*}
0=\langle [\nabla_X,\nabla_Y]Y-\nabla_{[X,Y]}Y,X\rangle
=\langle R^{\s^n}(X,Y)Y,X\rangle \neq 0
\end{align*}
as the right hand side is the sectional curvature of the plane spanned by \(X,Y\) which 
is non-vanishing due to 
$R^{\s^n}(X,Y)Z =\langle Y,Z\rangle X - \langle X,Z\rangle Y$, where $X,Y,Z$ are vector fields on $\s^n$, and the assumption that 
\(X\) is not a constant multiple of \(Y\). Hence, we arrive at a contradiction
such that there cannot exist a parallel vector field on \(\s^n,n\geq 2\).
\end{proof}

With these preparations at hand we are now able to examine the stability of the proper biharmonic map (\ref{dfq:biharmonic-q}).

\begin{Prop}
\label{q-stable}
Let \(m\geq 5.\)
The proper biharmonic map \(q\colon B^m\to \s^{m^\ell}\) defined as in
\eqref{dfq:biharmonic-q} is strictly stable if the following inequality is satisfied
\begin{align}
\label{eq:quadratic}
 m^2-12(\ell(\ell+m-2)+2m-8)>0.
\end{align}
\end{Prop}
\begin{proof}
Plugging the identities (\ref{prop-q}) into the stability condition
\eqref{eq:stability-intrinsic} we obtain
\begin{align*}
&\int_{B^m}\big(|\Delta\eta|^2+2|dq|^4|\eta|^2\big)\,dx-\int_{B^m}\langle\Delta^2q,q\rangle|\eta|^2\,dx 
-6\int_{B^m}|dq|^2|d\eta|^2\,dx\\   
&=\int_{B^m}|\Delta\eta|^2\,dx
-\delta\int_{B^m}\frac{|d\eta|^2}{r^2}\,dx,
\end{align*}
where 
\begin{align*}
\delta=3(\ell(\ell+m-2)+2m-8).
\end{align*}
Now, we employ \eqref{eq:hardy-laplace} and get 
\begin{align*}
&\int_{B^m}\big(|\Delta\eta|^2+2|dq|^4|\eta|^2\big)\,dx-\int_{B^m}\langle\Delta^2q,q\rangle|\eta|^2\,dx 
-6\int_{B^m}|dq|^2|d\eta|^2\,dx\\   
&\geq 
(\frac{m^2}{4}-\delta)\int_{B^m}\frac{|d\eta|^2}{r^2}dx.
\end{align*}
As spheres do not admit parallel vector fields, see Lemma\,\ref{lem:parallel_vf},  
we have $\int_{B^m}\frac{|d\eta|^2}{r^2}dx>0$. 
Therefore, and by inequality (\ref{eq:quadratic}), the map $q$ is strictly stable.
\end{proof}

Using Proposition\,\ref{q-stable} we can now prove the main result of this article.

\begin{proof}[Proof of Theorem \ref{thm:intrinsic}]
The statement now immediately follows from solving the quadratic equation \eqref{eq:quadratic}.    
\end{proof}

\bibliographystyle{plain}
\bibliography{mybib}

\end{document}